\RequirePackage{amsmath}
\documentclass[runningheads]{llncs}
\usepackage{graphicx}
\usepackage{color}
\usepackage{amssymb}
\usepackage{booktabs}
\usepackage{algorithm,algpseudocode}
\usepackage{hyperref}

\newcommand{\R}{\mathbb{R}}   
\newcommand{\Ruinf}{\mathbb{R} \cup \{ +\infty \}}
\newcommand{\clconv}{\overline{\text{conv}}}
\newcommand{\dom}{\operatorname{dom}}
\newcommand{\pfx}{\partial f(x)}
\newcommand{\eto}{\xi_{21}}
\newcommand{\ett}{\xi_{22}}
\newcommand{\etz}{\xi_{20}}
\newcommand{\eoo}{\xi_{11}}
\newcommand{\eot}{\xi_{12}}
\newcommand{\eoz}{\xi_{10}}
\newcommand{\ezo}{\xi_{01}}
\newcommand{\ezt}{\xi_{02}}
\newcommand{\ezz}{\xi_{00}}

\newcommand{\pgx}{\partial g(x)}
\newcommand{\prx}{\partial r(x)}

\newcommand{\ngx}{\nabla g(x)}
\newcommand{\nrx}{\nabla r(x)}

\newcommand{\intE}{\text{ri}(E)}

\newcommand{\intP}{\text{int}(P)}

\begin{document}
\title{Towards the biconjugate of bivariate piecewise quadratic functions\thanks{This work was supported in part by Discovery grants 298145-2013 and RGPIN-
		2018-03928 (second author) from the Natural Sciences and Engineering Research Council of Canada
		(NSERC) and by The University of British Columbia, Okanagan campus. Part of the research
		was performed in the Computer-Aided Convex Analysis (CA2) laboratory funded by a Leaders
		Opportunity Fund (LOF, John R. Evans Leaders Fund -- Funding for research infrastructure) from
		the Canada Foundation for Innovation (CFI) and by a British Columbia Knowledge Development
		Fund (BCKDF).}
		\thanks{This is a pre-copyedited version of a contribution published in
			Optimization of Complex Systems: Theory, Models, Algorithms and Applications edited by Hoai An Le Thi, Hoai Minh Le, Tao Pham Dinh, published by Springer Nature Switzerland AG. The definitive
			authenticated version is available online via \url{https://doi.org/10.1007/978-3-030-21803-4_27}}
}
%
%
\author{Deepak Kumar \and
Yves Lucet\inst{1}\orcidID{0000-0003-0349-6737}}
\authorrunning{D. Kumar et al.}
%
\institute{University of British Columbia Okanagan, 3187 University Way, Kelowna BC V1V 1V7 Canada
\email{yves.lucet@ubc.ca}, 
\url{https://people.ok.ubc.ca/ylucet/}}
\maketitle              

\begin{abstract}
Computing the closed convex envelope or biconjugate is the core operation that bridges the domain of nonconvex with convex analysis. We focus here on computing the conjugate of a bivariate piecewise quadratic function defined over a polytope. First, we compute the convex envelope of each piece, which is characterized by a polyhedral subdivision such that over each member of the subdivision, it has a rational form (square of a linear function over a linear function). Then we compute the conjugate of all such rational functions. It is observed that the conjugate has a parabolic subdivision such that over each member of its subdivision, it has a fractional form (linear function over square root of a linear function). This computation of the conjugate is performed with a worst-case linear time complexity algorithm.

Our results are an important step toward computing the conjugate of a piecewise quadratic function, and further in obtaining explicit formulas for the convex envelope of piecewise rational functions.

\end{abstract}
\section{Introduction}
Computational convex analysis (CCA) focuses on creating efficient algorithms to compute fundamental transforms arising in the field of convex analysis. Computing the convex envelope or biconjugate is the core operation that bridges the domain of nonconvex analysis with convex analysis. Development of most of the algorithms in CCA began with the Fast Legendre Transform (FLT) in~\cite{brenier1989algorithme}, which was further developed in~\cite{corrias1996fast,lucet1996fast}, and improved to the optimal linear worst-case time complexity in~\cite{lucet1997faster} and then~\cite{gardiner2011graph,lucet2006fast}. More complex operators were then considered~\cite{bauschke2008proximal,bauschke2008transform,lucet2009piecewise} (see~\cite{lucet2010shape} for a survey including a list of applications).

Piecewise Linear Quadratic (PLQ) functions (piecewise quadratic functions over a polyhedral partition) are well-known in the field of convex analysis~\cite{rockafellar2009variational} with the existence of linear time algorithms for various convex transforms~\cite{lucet2009piecewise,bauschke2008transform}. Computing the full graph of the convex hull of univariate  PLQ functions is possible in optimal linear worst-case time complexity~\cite{gardiner2010convex}.

For a function $f$ defined over a region $P$, the pointwise supremum of all its convex underestimators is called the \emph{convex envelope} and is denoted $	\clconv f_P(x,y)$. Computing the convex envelope of a multilinear function over a unit hypercube is NP-Hard~\cite{crama1989recognition}. However, the convex envelope of functions defined over a polytope $P$ and restricted by the vertices of $P$ can be computed in finite time using a linear program~\cite{tardella2004existence,tardella2008existence}. A method to reduce the computation of convex envelope of functions that are one lower dimension($\R^{n-1}$) convex and have indefinite Hessian to optimization problems in lower dimensions is discussed in~\cite{jach2008convex}.

Any general bivariate nonconvex quadratic function can be linearly transformed to the sum of bilinear and a linear function. Convex envelopes for bilinear functions over rectangles have been discussed in~\cite{mccormick1976computability} and validated in~\cite{al1983jointly}. The convex envelope over special polytopes (not containing edges with finite positive slope) was derived in~\cite{sherali1990explicit} while~\cite{linderoth2005simplicial} deals with bilinear functions over a triangle containing exactly one edge with finite positive slope. The convex envelope over general triangles and triangulation of the polytopes through doubly nonnegative matrices (both semidefinite and nonnegative) is presented in~\cite{anstreicher2012convex}.

In~\cite{locatelli2014technique}, it is shown that the analytical form of the convex envelope of some bivariate functions defined over polytopes can be computed by solving a continuously differentiable convex problem. In that case, the convex envelope is characterized by a polyhedral subdivision.

The Fenchel conjugate $f^*(s) = \sup_{x \in \R^n}[\langle s,x \rangle - f(x)]$ (we note $\langle s,x \rangle = s^T x$) of a function $f:\R^n \to \Ruinf$ is also known as the \emph{Legendre-Fenchel Transform} or \emph{convex conjugate} or simply \emph{conjugate}. It plays a significant role in \emph{duality} and computing it is a key step in solving the dual optimization problem~\cite{rockafellar2009variational}. Most notably, the biconjugate is also the closed convex envelope. 

A method to compute the conjugate known as the fast Legendre transform was introduced in~\cite{brenier1989algorithme} and studied in~\cite{corrias1996fast,lucet1996fast}. A linear time algorithm was later introduced by Lucet to compute the discrete Legendre transform~\cite{lucet1997faster}. Those algorithms are numeric and do not provide symbolic expressions.

Computation of the conjugate of convex univariate PLQ functions have been well studied in the literature and linear time algorithms have been developed in~\cite{gardiner2013computing,gardiner2014computing}. Recently, a linear time algorithm to compute the conjugate of convex bivariate PLQ functions was proposed in~\cite{haque2018linear}.

Let $f:\R^n \to \Ruinf$ be a piecewise function, i.e. $f(x)=f_i(x)$ if $x\in P_i$ for $i=1,\dots,N$. From~\cite[Theorem 2.4.1]{hiriart1993convex}, we have
$(\inf_i f_i)^* = \sup_i f^*_i$, and from~\cite[Proposition 2.6.1]{hiriart1993convex},
$	\clconv (\inf_i (f_i + I_{P_i})) = \clconv ( \ \inf_i [\clconv (f_i + I_{P_i}) ] \ ) $where $I_{P_i}$ is the indicator function for $P_i$.
Hence, 
$\clconv (\inf_i (f_i + I_{P_i})) = ( \ \sup_i \ [\clconv(f_i + I_{P_i})]^* \ )^* $. This provides an algorithm to compute the closed convex envelope: (1) compute the convex envelope of each piece, (2) compute the conjugate of the convex envelope of each piece, (3) compute the maximum of all the conjugates, and (4) compute the conjugate of the function obtained in (3) to obtain the biconjugate. The present work focuses on Step (2).

Recall that given a quadratic function over a polytope, the eigenvalues of its symmetric matrix determine how difficult its convex envelope is to compute (for computational purposes, we can ignore the affine part of the function). If the matrix is semi-definite (positive or negative), the convex envelope is easily computed. When it is indefinite, a change of coordinate reduces the problem to finding the convex envelope of the function $(x,y)\mapsto xy$ over a polytope, for which step (1) is known~\cite{locatelli2016polyhedral}.

The paper is organized as follow. Section~\ref{s:dom} focuses on the domain of the conjugate while Section~\ref{s:conj} determines the symbolic expressions. Section~\ref{s:conc} concludes the paper with future work.

\section{Preliminaries and notations}\label{s:pre}
The \emph{subdifferential} $\pfx$ of a function $f:\R^n \to \Ruinf$ at any $x \in \dom(f)=\{x : f(x) < \infty\}$ is
	$\pfx =  \{ s : f(y) \ge f(x) + \langle s, y-x \rangle , \forall y \in \dom(f) \}$ ($\pfx = \{ \nabla f(x) \}$ when $f$ is differentiable at $x$). We note $I_P$ the indicator function of the set $P$, i.e. $I_P(x)=0$ when $x\in P$ and $I_P(x)=+\infty$ when $x\notin P$.

	A \emph{parabola} is a two dimensional planar curve whose points $(x,y)$ satisfy the equation $a x^2 + bxy + cy^2 + dx + ey + f = 0$ with $b^2 - 4ac = 0$. A \emph{parabolic region} is formed by the intersection of a finite number of parabolic inequalities, i.e. $P_r = \{ x \in \R^2 : C_p^i(x) \le 0, i \in \{1,\cdots,k \}\}$ where $C_p^i(x)=a_i x_1^2 + b_ix_1x_2 + c_ix_2^2 + d_ix_1 + e_ix_2 + f_i$ and $b_i^2 - 4 a_i c_i = 0$. The set $P_r^i = \{ x \in \R^2 : C_p^i(x) \le 0 \}$ is convex, but $P_s^i = \{  x \in \R^2 : C_p^i(x) \ge 0 \}$ is not.

	A convex set $R = \bigcup_{i \in \{1,...,m \}} R_i,$ $ R \subseteq \R^2 $, defined as the union of a finite number of parabolic regions is said to have a \emph{parabolic subdivision} if for any $j,k \in \{1,\cdots,m\}, \ j \ne k, R_j \bigcap R_k$ is either empty or is contained in a parabola.

\section{The domain of the conjugate}\label{s:dom}
Given a nonconvex PLQ function, we first compute the closed convex envelope of each piece and obtain a piecewise rational function~\cite{locatelli2016polyhedral}. We now compute the conjugate of such a rational function over a polytope by first computing its domain, which will turn out to be a parabolic subdivision. Recall that for PLQ functions, $\dom f^* = \partial f (\dom f)$. We decompose the polytope $\dom f=P $ into its interior, its vertexes, and its edges. 

Following~\cite{locatelli2016polyhedral}, we write a rational function as
\begin{equation}
r(x,y) =  \frac{(\xi_1 (x,y))^2}{\xi_2(x,y)} + \xi_0(x,y), 
\label{eq:r}
\end{equation}
where $\xi_i(x,y)$ are linear functions in $x$ and $y$.

\begin{proposition}[Interior]
	\label{pro:paraCurve}
	Consider $r$ defined by \eqref{eq:r}, there exists $ \alpha_{ij} $ such that $ \bigcup_{ x \in \dom(r)} \partial r(x) = \{ s : C_r(s) = 0 \}$, where $C_r(s) = \alpha_{11}s_1^2 + \alpha_{12} s_1 s_2 + \alpha_{22} s_2^2 + \alpha_{10} s_1 + \alpha_{02} s_2 + \alpha_{00}$ and $\{s:C_r(s)=0\}$ is a parabolic curve.
\end{proposition}
\begin{proof}
	Note $ \xi_1(x) = \eoo x_1 + \eot x_2 + \eoz, \ \xi_2(x) =  \eto x_1 + \ett x_2 + \etz$ and $ \xi_0(x) =  \ezo x_1 + \ezt x_2 + \ezz $. Since $r$ is differentiable everywhere in $\dom(r) =  \R^2 / \{ z:  \xi_2(z) = 0 \},$ for any $x \in \dom(r)$ we compute $s = \nabla r(x)$ as $s_i = 2 \xi_{1i} t - \xi_{2i} t^2 + \xi_{0i}$ for $i=1,2$, where $t = ({\eoo x_1 + \eot x_2 + \eoz})/({\eto x_1 + \ett x_2 + \etz})$. Hence, $s=\nabla r(x)$ represents the parametric equation of a conic section, and by eliminating $t$, we get
$C_r(s)=0$ where
	\[
C_r(s) =\alpha_{11}s_1^2 + \alpha_{12} s_1 s_2 + \alpha_{22} s_2^2 + \alpha_{10} s_1 + \alpha_{02} s_2 + \alpha_{00},
\]
with $\alpha_{11} = \eto^2 \ett^2, \ \alpha_{12} = -2 \eto^3 \ett, \ \alpha_{22} = \eto^4 $ and other $\alpha_{ij}$ are functions of the coefficients of $r$. We check that $		\alpha_{12}^2 - 4 \alpha_{11} \alpha_{22} = (-2 \eto^3 \ett)^2 - 4 \eto^6 \ett^2 = 0$, so
	the conic section is a parabola. Consequently, for all $x \in \dom(r), \partial r(x)$ is contained in the parabolic curve $C_r(s)=0$, i.e.
\[
	 	\bigcup_{ x \in \dom(r)} \partial r(x) \subset \{ s: C_r(s) = 0  \}.
\]	

	Conversely, any point $s_r$ that satisfies $C_r(s_r) = 0$, satisfies the parametric equation as well, so the converse inclusion is true.
\end{proof}

\begin{corollary}[Interior]
	\label{cor:rc}
		For a bivariate rational function $r$, and a polytope $P$, define $f(x) = r(x) + I_P(x)$, then for all $ x \in \text{int}(P)$, the set $ \bigcup_{x \in \text{int}(P)}$ $\partial f(x)$ is contained inside a parabolic arc.
\end{corollary}
\begin{proof}
	 We have, $\bigcup_{x \in \text{int}(P)} \partial f(x) \subseteq \bigcup_{x \in \intP} \partial r(x)$ and $\bigcup_{x \in \intP} \prx$ $ \subset \mathcal{P}$ where $\mathcal{P} \subset \R^2$ is a parabolic curve (from Proposition~\ref{pro:paraCurve}). Since $P$ is connected, we obtain that $\bigcup_{x \in \intP} \prx$ is contained in a parabolic arc. 
\end{proof}

Next we compute the subdifferential at any vertex in the smooth case (the proof involves a straightforward computation of the normal cone).
\begin{lemma}[Vertices]
	\label{lemm:verNc}
	For $g \in \mathcal{C}^1$, $P$ a polytope, and $v$ vertex. Let $f(x) = g(x) + I_P(x)$. Then $\partial f(v)$ is an unbounded polyhedral set.
\end{lemma}

There is one vertex at which both numerator and denominator equal zero although the rational function can be extended by continuity over the polytope; we conjecture the result based on numerous observations.
\begin{conjecture}[Vertex]
	\label{conj:undfPara}
	Let $r$ as in \eqref{eq:r}, $f(x) = r(x)+I_P(x)$ and $v$ be a vertex of $P$ with $\xi_1(v) = \xi_2(v)=0$. Then $\partial f(v)$ is a parabolic region.  
\end{conjecture}

\begin{lemma}[Edges]
	For $g \in \mathcal{C}^1$, a polytope $P$, and an edge $E = \{ x: x_2 = m x_1 + c,  x_{l_1} \le x_1 \le x_{u_1}  \}$ between vertices $x_l$ and $x_u$, let $f(x) = g(x) + I_P(x)$, then $\bigcup_{x \in \intE} \partial f(x) = \bigcup_{x \in \intE} \{ s + \ngx: s_1 + m s_2 = 0,  s_2 \ge 0 \} $.
	\label{lemm:sdedge}
\end{lemma}
\begin{proof}
	For all $x \in \intE$, $\pfx = \pgx + N_P(x).$ Let $L(x) = x_2 - m x_1 - c$ be the expression of the line joining $x_l$ and $x_u$ such that $P \subset \{ x : L(x) \le 0  \} $. (The case $ P \subset \{ x: L(x) \ge 0 \}$ is analogous. )
	
	Since $P \subset \R^2$ is a polytope, for all $x \in \intE, \ N_P(x) = \{s : s = \lambda \nabla L(x), \ \lambda \ge 0 \}$ is the normal cone of $P$ at $x$ and can be written $N_P(x) = \{ s: s_1 + m s_2 = 0, \ s_2 \ge 0 \}$.
	In the special case when $E = \{ x:  x_1 = d, x_{l_1} \le x_1 \le x_{u_1} \}$, $L(x) = x_1 -d $ and $N_P(x) = \{ s: s_2 = 0, s_1 \ge 0 \}$.
	Now for any $ x \in \intE$, $\pfx = \pgx + N_P(x) = \{ s + \ngx : s_1 + m s_2 = 0, \ s_2 \ge 0 \}$, so
	\begin{equation*}
		\begin{aligned}
		\bigcup_{x \in \intE } \pfx &= \bigcup_{x \in \intE } \{ s + \ngx : s_1 + m s_2 = 0, \ s_2 \ge 0 \}. 
		\end{aligned}
	\end{equation*}
\end{proof}

\begin{proposition}[Edges]
	\label{pro:rdedge}
	For $r$ as in \eqref{eq:r}, a polytope $P$ and an edge $E = \{ x: x_2 = m x_1 + c , v_1^- \le x_1 \le v_1^+ \}$ between vertices $v^-$ and $v^+$, let $f(x) = r(x) + I_P(x)$, then $\bigcup_{x \in \intE} \pfx$ is either a parabolic region or a ray.
\end{proposition}
\begin{proof}
	From Corollary~\ref{cor:rc}, there exists $l,u \in \R^2$ such that $\bigcup_{x \in \intE } \prx = \bigcup_{x \in \intE } \{ s: C_r(s) = 0, l_1 \le s_1 \le u_1 \}$. So computing $\bigcup_{x \in \intE } \pfx$ leads to the following two cases:
	
	\textbf{Case 1 ($ l = u$)} Same case as when $r$ is quadratic (known result).
	
	\textbf{Case 2 ($l \ne u$)} By setting $g=r$ in Lemma~\ref{lemm:sdedge}, for any $x \in \intE,$ $\pfx = \{ s : s_1 + m s_2 = 0, s_2 \ge 0 \}$. Similar to the quadratic case, when $\nrx = l, \ \pfx = \{ s  : s_1 + m s_2 -(l_1 + m l_2) = 0, s_2 \ge l_2  \}$ and when $\nrx = u, \ \pfx = \{ s  : s_1 + m s_2 -(u_1 + m u_2) = 0, s_2 \ge u_2  \}$. Assume $\pfx \subset \{ s: C_r(s) \le 0 \}$ (the case $\pfx \subset \{ s: C_r(s) \ge 0  \}$ is analogous). Then
	\begin{equation*}
		\begin{aligned}
			\bigcup_{x \in \intE } \pfx &= \bigcup_{x \in \intE } \{ s + \nabla r(x) : s_1 + m s_2 = 0, s_2 \ge 0 \} \\
			&= \{ s: l_1 + m l_2 \le s_1 + m s_2 \le u_1 + m u_2, \ C_r(s) \le 0 \}
		\end{aligned}
	\end{equation*}
	is a parabolic region.
\end{proof}

By gathering Lemma~\ref{lemm:verNc}, Proposition~\ref{pro:rdedge}, and Corollary~\ref{cor:rc}, we obtain.

\begin{theorem}	[Parabolic domain]
	\label{theo:rparadom}
	Assuming Conjecture~\ref{conj:undfPara} holds, $r$ is as in \eqref{eq:r}, $P$ is a polytope, and $f(x) = r(x) + I_P(x)$. Then $\bigcup_{x \in P} \partial f(x)$ has a parabolic subdivision.
\end{theorem}

\begin{figure}[t]
	\includegraphics[width=\linewidth]{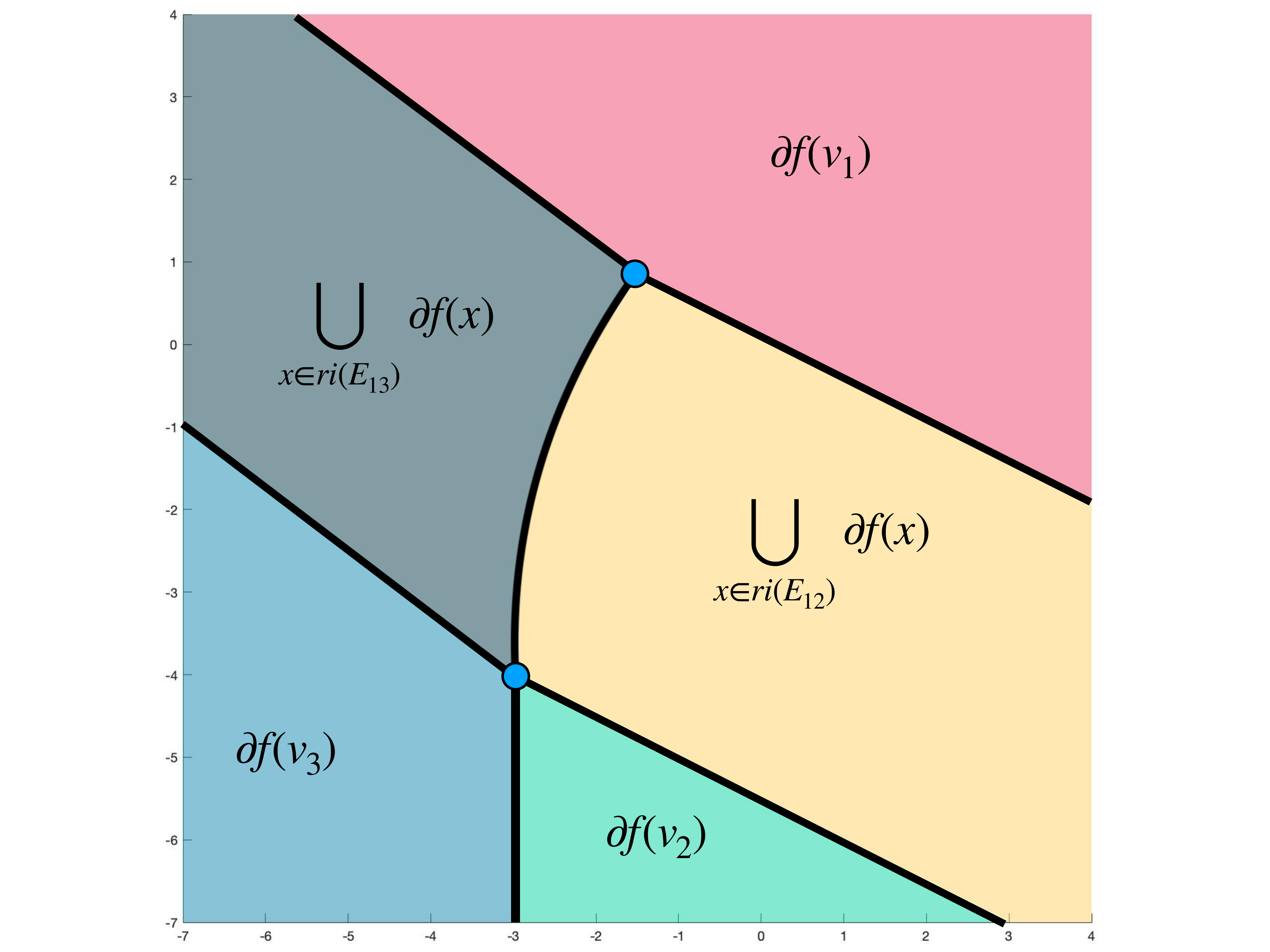}
	\caption{Parabolic subdivision for $r$ and $P$ from Example~\ref{fig:ex2rP1}}
	\label{fig:ex2rP1}
\end{figure}

\begin{example}
	\label{ex:rfP1}
	For $ \displaystyle r = \frac{36x_1^2 + 21 x_1 x_2 + 36 x_2^2 - 81 x_1 + 24 x_2 - 252}{-12 x_1 + 9 x_2 + 75} $ and polytope $P$ formed by vertices $v_1 = (-1,1), v_2 = (-3,-3)$ and $v_3 = $ $(-4,-3)$, let $f(x) = r(x) + I_P(x)$. We have $ \bigcup_{ x \in \dom(r)} \prx$ $= \{ s: C(s)=0 \}$ where $C(s) = 9s_1^2 + 24s_1s_2 - 234s_1 + 16s_2^2 + 200s_2 - 527$. 
	%
	The parabolic subdivision for this example is shown in Figure~\ref{fig:ex2rP1}.
\end{example}

\section{Conjugate Expressions}\label{s:conj}
Now that we know $\dom f^*$ as a parabolic subdivision, we turn to the computation of its expression on each piece. We note
\begin{align}
	g_f(s_1,s_2) &= \frac{\psi_1(s_1,s_2)}{\zeta_{00} \sqrt{\psi_{1/2}(s_1,s_2)}} + \psi_0(s_1,s_2) \label{eq:gff} \\
	g_q(s_1,s_2) &= \zeta_{11} s_1^2 + \zeta_{12}s_1s_2 + \zeta_{22}s_2^2 +  \zeta_{10} s_1 + \zeta_{01} s_2 + \zeta_{00} \label{eq:gfq} \\
	g_l(s_1,s_2) &= \zeta_{10} s_1 + \zeta_{01} s_2 + \zeta_{00} \label{eq:gfl}
\end{align}
where $\psi_0,\psi_{1/2}$ and $\psi_1$ are linear functions in $s$, and $\zeta_{ij} \in \R$.

\begin{theorem}
	Assume Conjecture~\ref{conj:undfPara} holds. For $r$ as in \eqref{eq:r}, a polytope $P$, and $f(x) = r(x) + I_P(x)$, the conjugate $f^*(s)$ has a parabolic subdivision such that over each member of its subdivision it has one of the forms in~\eqref{eq:gff}-\eqref{eq:gfl}
	\label{theo:rconjxps}
\end{theorem}
\begin{proof}
	We compute the critical points for the optimization problem defining $f^*$.
	
	\textbf{Case 1 (Vertices)} For any vertex $v$, $f^*(s) = s_1 v_1 + s_2 v_2 - r(v)$
	is a linear function of form~\eqref{eq:gfl} defined over an unbounded polyhedral set (from Lemma~\ref{lemm:verNc}). In the special case, when $\partial f(v)$ is a parabolic region (Conjecture~\ref{conj:undfPara}), the conjugate would again be a linear function but defined over a parabolic region.
	
	\textbf{Case 2 (Edges)} Let $F$ be the set of all the edges, and $E = \{ x : x_2 = m x_1 + c, l_1 \le x_1 \le u_1 \} \in F$ be an edge between vertices $l$ and $u$, then $	f^*(s) = \sup_{x \in \intE} \{ \langle s,x \rangle - (r(x) + I_P(x))  \}$.
	By computing the critical points, we have $ s -(\nabla r(x) + N_P(x)) = 0$ where $N_P(x) = \{ s: s= \lambda (-m,1), \lambda \ge 0 \}$ with $m$ the slope of the edge. So
	\begin{equation}
		\begin{aligned}
			s_1 &= - \eto t^2 + 2 \eoo t  + \ezo - m\lambda \\
			s_2 &= - \ett t^2 + 2 \eot t  + \ezt + \lambda
		\end{aligned}
		\label{eq:rxp1}
	\end{equation}
	where $\displaystyle t = \frac{\eoo x_1 + \eot x_2 + \eoz}{\eto x_1 + \ett x_2 + \etz}$. Since $x \in \intE$, we have 
\begin{equation}
	x_2 = m x_1 + c,
	\label{eq:rxp2}
\end{equation}
	which with~\eqref{eq:rxp1} gives
	\begin{equation}
	\begin{aligned}
	x_1 &= \begin{cases}
	\gamma_{10} s_1 + \gamma_{01} s_2 + \gamma_{00} & \text{ when } \eto + m \ett =0 \\
	\displaystyle 
	\frac{\gamma_{00} \pm \gamma_{1/2} \sqrt{\gamma_{10/2} s_1 + \gamma_{01/2} s_2 + \gamma_{00/2}}}{\pm \gamma_{-1/2} \sqrt{\gamma_{10/2} s_1 + \gamma_{01/2} s_2 + \gamma_{00/2}}} & \text{ otherwise, }
	\end{cases}
	\end{aligned}
	\label{eq:rxpx1}
	\end{equation}
	where all $\gamma_{ij}$ and $\gamma_{ij/k}$ are defined in the coefficients of $r$, and parameters $m$ and $c$. When $\eto + m \ett \ne 0$, solving~\eqref{eq:rxp1} and \eqref{eq:rxp2}, leads to a quadratic equation in $t$ with coefficients as linear functions in $s$.

	By substituting~\eqref{eq:rxpx1} and~\eqref{eq:rxp2} in $f^*(s)$, when $\eto + m \ett \ne 0$, we have
	$$f^*(s) = \frac{\psi_1(s_1,s_2)}{\zeta_{00} \sqrt{\psi_{1/2}(s_1,s_2)}} + \psi_0(s_1,s_2),$$
	and when $\eto + m \ett = 0$,
	$$f^*(s) = \zeta_{11} s_1^2 + \zeta_{12}s_1s_2 + \zeta_{22}s_2^2 +  \zeta_{10} s_1 + \zeta_{01} s_2 + \zeta_{00}, $$
	where all $\zeta_{ij},\psi_{i}$ and $\psi_{i/j}$ are defined in the coefficients of $r$, and parameters $m$ and $c$, with $\psi_i(s)$ and $\psi_{i/j}(s)$ linear functions in $s$. 
	
	 From Proposition~\ref{pro:rdedge}, $\bigcup_{x \in \intE} \pfx$ is either a parabolic region or a ray. So for any $E$, the conjugate is a fractional function of form~\eqref{eq:gff} defined over a parabolic region. When $\bigcup_{ x \in \intE} \pfx$ is a ray, the computation of the conjugate is deduced from its neighbours by continuity.
	
	\textbf{Case 3 (Interior)} Since $\bigcup_{x \in \intP } \pfx$ is contained in a parabolic arc (from Corollary~\ref{cor:rc}), the computation of the conjugate is deduced by continuity. 
\end{proof}

\begin{example}\label{sec:exratz}
For a bivariate rational function $\displaystyle r(x) = \frac{x_2^2}{x_2 - x_1 + 1}$ defined over a polytope $P$ with vertices $v_1 = (1,1), v_2 = (1,0)$ and $v_3 = (0,0)$, let $f(x) = r(x) + I_P(x)$. 

	\begin{figure}
		\includegraphics[width=\linewidth]{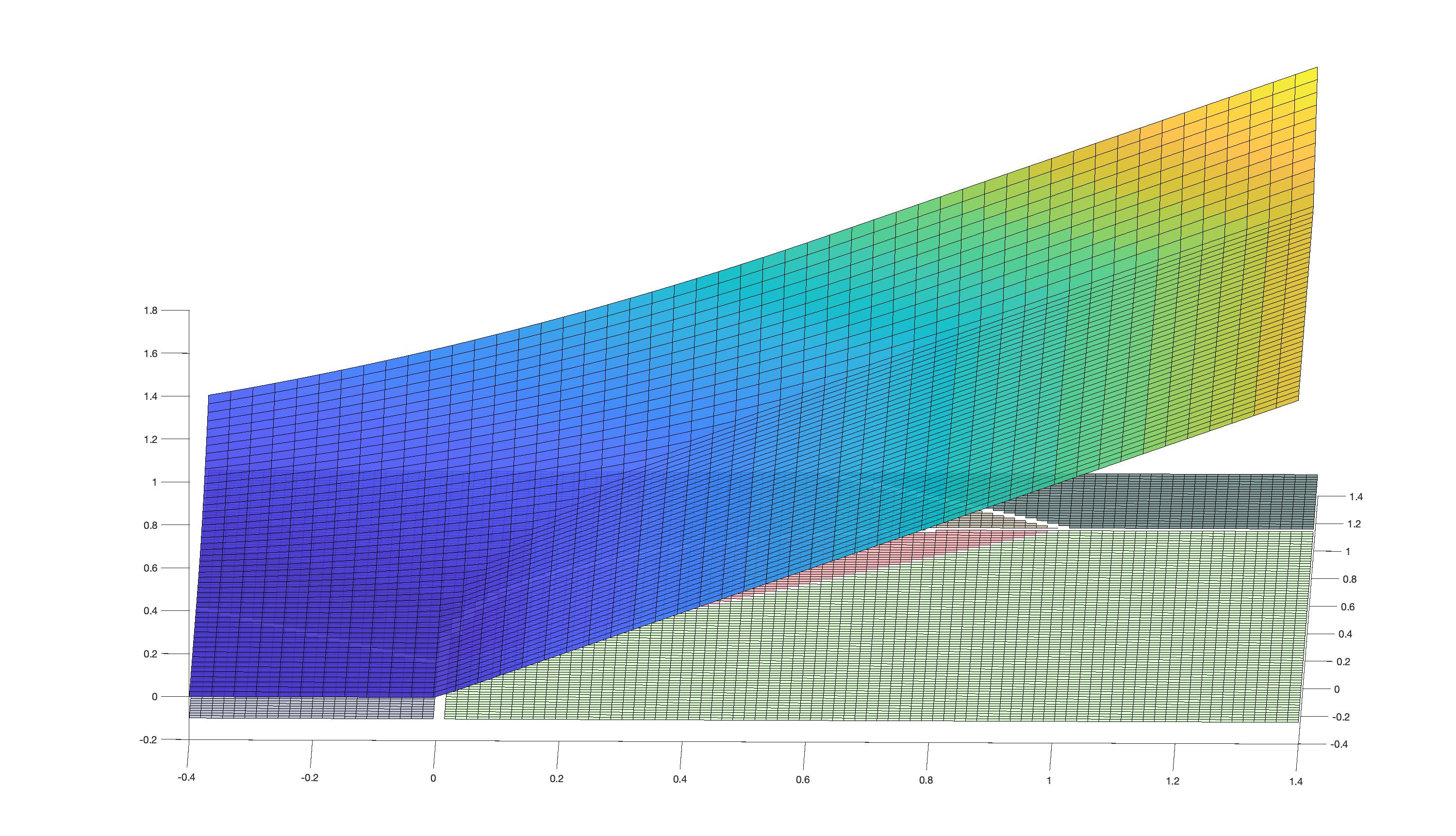}
		\caption{Conjugate for Example~\ref{sec:exratz}}
		\label{fig:conjEx1rz}
	\end{figure}

	The conjugate (shown in Figure~\ref{fig:conjEx1rz}) can be written 
	\begin{align*}
		f^*_P(s) = \begin{cases}
		0 & s \in R_1 \\
		s_1 & s \in  R_2 \\
		s_1 & s \in R_3 \\
		s_1 + s_2 -1 & s \in R_4 \\
		\frac{1}{4}(s_1 + s_2)^2 & s \in R_5
		\end{cases}
	\end{align*}
	where
	\begin{align*}
		R_1 &= \{ s: s_2 \ge -s_1 + 2, s_2 \ge  1 \} \\
		R_2 &= \{ s: s_2 \ge s_1, s_1^2 + 2s_1s_2 - 4s_1 + s_2^2 \le 0 \} \\
		R_3 &= \{ s: s_2 \le s_1, s_2 \le 1,  s_1 \ge 0 \} \\
		R_4 &= \{ s: 0 \le s_1, s_2 \le -s_1 \} \\
		R_5 &= \{ s: s_2 \ge -s_1, s_2 \le -s_1 + 2, s_2 \ge s_1, s_1^2 + 2s_1s_2 - 4s_1 + s_2^2 \ge 0 \}.
	\end{align*}
\end{example}
	
\section{Conclusion and future work}\label{s:conc}
\begin{figure}
		\includegraphics[width=\linewidth]{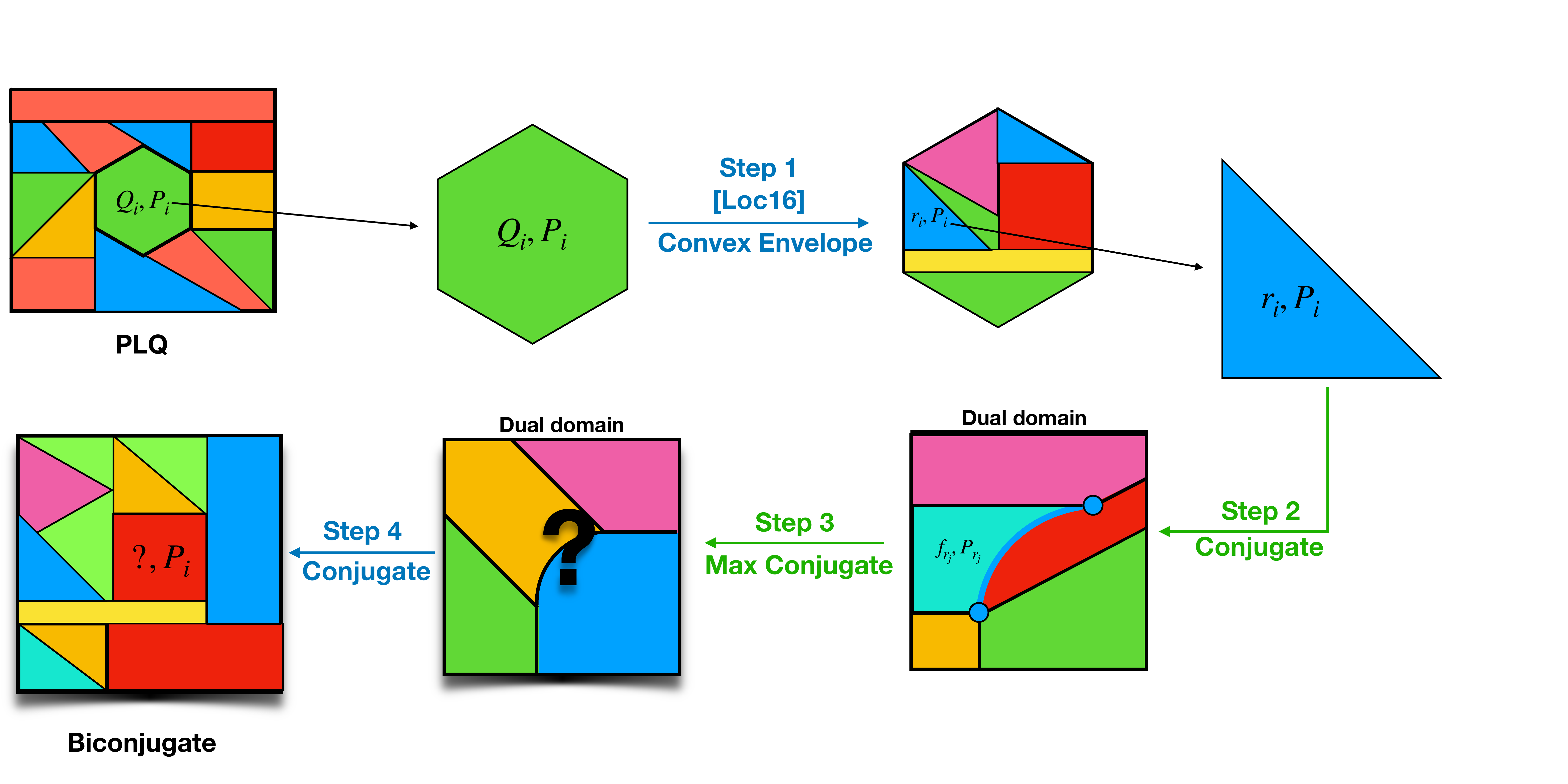}
		\caption{Summary}
		\label{fig:summary}
\end{figure}

Figure~\ref{fig:summary} summarizes the strategy. Given a PLQ function, for each piece, its convex envelope is computed as the convex envelope of a quadratic function over a polytope using \cite{locatelli2014technique}. This is the most time consuming operation since the known algorithms are at least exponential. For each piece, we obtain a piecewise rational function. Then we take each of those pieces, and compute its conjugate to obtain a fractional function over a parabolic subdivision. That computation is complete except for Conjecture~\ref{conj:undfPara}. Note that there is only a single problematic vertex $v$ and since the conjugate is full domain, we can deduce $\partial f(v)$ by elimination.

Future work will focus on Step 3, which will give the conjugate of the original PLQ function. This will involve solving repeatedly the map overlay problem and is likely to take exponential time. From hundred of examples we ran, we expect the result to be a fractional function of unknown kind over a parabolic subdivision; see Figure~\ref{fig:summary}, bottom row, middle figure. The final step will be to compute the biconjugate (bottom-left in Figure~\ref{fig:summary}. We know it is a piecewise function over a polyhedral subdivision but do not know the formulas.

 \bibliographystyle{splncs04}
 \bibliography{lucet}
%
%
%
%
%
\end{document}